\documentclass[authoryear]{elsarticle}

\usepackage{amsmath}
\usepackage{amsfonts}
\usepackage{bm}
\usepackage{tikz}
\usepackage{pgfplots}
\usepackage{color, colortbl}
\usepackage{graphicx}

\usepackage{hyperref}

\newtheorem{thm}{Theorem}

\newproof{pf}{Proof}

\journal{arXiv}

\begin{document}

\begin{frontmatter}

\title{Numerical solution of time-dependent problems with fractional power elliptic operator}

\author[nsi,univ]{Petr N. Vabishchevich\corref{cor}}
\ead{vabishchevich@gmail.com}

\address[nsi]{Nuclear Safety Institute, Russian Academy of Sciences, 52, B. Tulskaya, Moscow, Russia}
\address[univ]{Peoples' Friendship University of Russia (RUDN University),
        		6 Miklukho-Maklaya St., 117198 Moscow, Russia}

\cortext[cor]{Corresponding author}

\begin{abstract}
An unsteady problem is considered for a space-fractional equation in a bounded domain. 
A first-order evolutionary equation involves a fractional power
of an elliptic operator of second order. Finite element approximation in
space is employed. To construct approximation in time, standard two-level schemes are used.
The approximate solution at a new time-level is obtained as a solution of
a discrete problem with the fractional power of the elliptic operator.
A Pade-type approximation is constructed on the basis of special
quadrature formulas for an integral representation of the fractional power elliptic operator
using explicit schemes. A similar approach is applied in the numerical implementation of implicit schemes. 
The results of numerical experiments are presented for a test two-dimensional problem.
\end{abstract}


\end{frontmatter}

\section{Introduction}

Many applied mathematical models involve both sub-diffusion (fractional in time) and super-diffusion 
(fractional in space) operators (see, e.g., \cite{podlubny1998fractional,uchaikin}). 
Super-diffusion problems are treated as evolutionary problems with a fractional power of an elliptic operator.
For example, suppose that in a bounded domain $\Omega$ on the set of functions 
$u(\bm x) = 0, \ \bm x \in \partial \Omega$, 
there is defined the operator $\mathcal{A}$: $\mathcal{A}  u = - \triangle u, \ \bm x \in \Omega$. 
We seek the solution of the Cauchy problem for the equation with the fractional power elliptic operator:
\[
 \frac{d u}{d t} + \mathcal{A}^{\alpha } u = f(t),
 \quad 0 < \alpha < 1, 
 \quad 0 < t \leq T, 
\] 
\[
 u(0) = u_0,
\] 
for a given $f(\bm x, t)$, $u_0(\bm x), \ \bm x \in \Omega$ using the notation $f(t) = f(\cdot,t)$. 

For approximation in space, we can apply finite volume or finite element methods
oriented to using arbitrary domains and irregular computational grids (\cite{KnabnerAngermann2003,QuarteroniValli1994}).
After this, we formulate the corresponding Cauchy problem with a fractional power of a
self-adjoint positive definite discrete elliptic operator 
(see \cite{bonito2015numerical,szekeres2016finite}) ---  a fractional power of a
symmetric positive definite matrix (\cite{higham2008functions}).

In the study of difference schemes for time-dependent problems of BVP for PDE, the general theory of
stability (well-posedness) for operator-difference schemes (\cite{Samarskii1989,SamarskiiMatusVabischevich2002})
is in common use. At the present time,
the exact (matching necessary and sufficient) conditions for stability are obtained for
a wide class of two- and three-level difference schemes considered in
finite-dimensional Hilbert spaces. We emphasize a constructive nature of the general theory of stability
for operator-difference schemes, where stability criteria are
formulated in the form of operator inequalities, which are easy to verify.
In particular, the backward Euler scheme and Crank-Nicolson 
scheme are unconditionally stable for a non-negative operator. 

Problems in numerical solving unsteady problems with fractional powers of operators appear
in using the simplest explicit approximations in time.
A practical implementation of such approach requires the matrix function-vector multiplication.
For such problems, different approaches (see \cite{higham2008functions}) are available.
Algorithms for solving systems of linear equations associated with
fractional elliptic equations that are based on Krylov subspace methods with the Lanczos approximation
are discussed, e.g., in \cite{ilic2009numerical}.
A comparative analysis of the contour integral method, the extended Krylov subspace method, and the preassigned 
poles and interpolation nodes method for solving
space-fractional reaction-diffusion equations is presented in \cite{burrage2012efficient}.
The simplest variant is associated with the explicit construction of the solution using the 
eigenvalues and eigenfunctions of the elliptic operator with diagonalization of the corresponding matrix
(\cite{bueno2012fourier,ilic2006numerical}). 
Unfortunately, all these approaches demonstrate very high computational complexity for multidimensional problems.

There does exist a general approach to solve approximately equations involving a fractional power of operators
based on an approximation of the original operator and then taking fractional power of its discrete variant. 
Using the Dunford-Cauchy formula the elliptic operator is represented as a contour integral  in the complex plane. 
Further, applying appropriate quadratures with integration nodes in the complex plane, it is necessary to select 
a proper method that involves only inversion of the original operator. 
The approximate operator is treated as a sum of resolvents (\cite{gavrilyuk2004data,gavrilyuk2005data})
ensuring the exponential convergence of quadrature approximations.
In \cite{bonito2015numerical}, there was presented a more promising variant of using quadrature
formulas with nodes on the real axis, which are constructed
on the basis of the corresponding integral representation 
for the power operator (see \cite{krasnoselskii1976integral,carracedo2001theory}).
In this case, the inverse operator of the problem has an additive representation, where
each term is an  inverse of the original elliptic operator.
A similar rational approximation to the fractional Laplacian operator 
is studied in \cite{AcetoNovat}.

We have proposed (\cite{vabishchevich2014numerical}) a numerical algorithm to solve
an equation for fractional power elliptic operators that is based on
a transition to a pseudo-parabolic equation.
For an auxiliary Cauchy problem, the standard two-level schemes are applied.
The computational algorithm is simple for practical use, robust, and applicable to solving
a wide class of problems. A small number of time steps is required to find a solution.
This computational algorithm for solving equations with fractional power operators
is promising for transient problems. 

As for one-dimensional problems for the space-fractional diffusion equation,
an analysis of stability and convergence for this equation was conducted in \cite{jin2014error} using
finite element approximation in space. A similar study for the Crank-Nicolson scheme was conducted earlier
in \cite{tadjeran2006second} using finite difference approximations in space.  
We highlight separately the works \cite{huang2008finite,sousa2012second,meerschaert2004finite}, 
where numerical methods for solving one-dimensional transient problems of convection and 
space-fractional diffusion equation are considered. 

In \cite{vabishchevich2016numerical}, an unsteady problem is considered for a space-fractional diffusion equation 
in a bounded domain.
To construct approximation in time, regularized two-level schemes are used (see \cite{Vabishchevich2014}).
The numerical implementation is based on solving the equation with 
the fractional power of the elliptic operator using an auxiliary Cauchy problem for a pseudo-parabolic equation
(\cite{vabishchevich2014numerical}).
Some more general unsteady problems are considered in
\cite{vabishchevich2016cnumerical,vabishchevich2016bnumerical}. 

In the present work, standard two-level schemes are applied to solve numerically a Cauchy problem for an evolutionary
equation of first order with a fractional power of the operator.
The numerical implementation is based on the rational approximation of the operator at a new time-level.
When implementing the explicit scheme, the fractional power of the operator is approximated on the basis of 
Gauss-Jacobi quadrature formulas for the corresponding integral representation.
In this case, we have (see \cite{frommer2014efficient}) a Pade-type approximation of the power function with 
a fractional exponent. A similar approach is used when considering implicit schemes.

The paper is organized as follows.
The formulation of an unsteady problem containing a fractional power
of an elliptic operator is given in Section 2. Here 
finite element approximations in space are also discussed. 
In Section 3, we construct the explicit  approximation in time
and investigate its stability. 
The numerical implementation is based on the rational approximation of
the fractional power operator. Implicit schemes are considered in Section 4.
The results of numerical experiments are described in Section 5.
At the end of the work the main results of our study are summarized.

\section{Problem formulation}

In a bounded polygonal domain $\Omega \subset \mathbb{R}^d$, $d=2,3$ with the Lipschitz continuous boundary $\partial\Omega$,
we search the solution for the problem with a fractional power of an elliptic operator.
Define the elliptic operator as
\begin{equation}\label{1}
  \mathcal{A}  u = - {\rm div}  k({\bm x}) {\rm grad} \, u + c({\bm x}) u
\end{equation} 
with coefficients $0 < \underline{k} \leq k({\bm x}) \leq \overline{k}$, $c({\bm x}) \geq 0$.
The operator $\mathcal{A}$ is defined on the set of functions $u({\bm x})$ that satisfy
on the boundary $\partial\Omega$ the following conditions:
\begin{equation}\label{2}
  k({\bm x}) \frac{\partial u }{\partial n } + g ({\bm x}) u = 0,
  \quad {\bm x} \in \partial \Omega ,
\end{equation} 
where $g ({\bm x}) \geq 0, \  {\bm x} \in \partial \Omega$.

In the Hilbert space $H = L_2(\Omega)$, we define the scalar product and norm in the standard way:
\[
  (u,v) = \int_{\Omega} u({\bm x}) v({\bm x}) d{\bm x},
  \quad \|u\| = (u,u)^{1/2} .
\] 
For the spectral problem
\[
 \mathcal{A}  \varphi_k = \lambda_k \varphi_k, 
 \quad \bm x \in \Omega , 
\] 
\[
  k({\bm x}) \frac{\partial  \varphi_k}{\partial n } + g ({\bm x}) \varphi_k = 0,
  \quad {\bm x} \in \partial \Omega , 
\] 
we have 
\[
 \lambda_1 \leq \lambda_2 \leq ... ,
\] 
and the eigenfunctions  $ \varphi_k, \ \|\varphi_k\| = 1, \ k = 1,2, ...  $ form a basis in $L_2(\Omega)$. Therefore, 
\[
 u = \sum_{k=1}^{\infty} (u,\varphi_k) \varphi_k .
\] 
Let the operator $\mathcal{A}$ be defined in the following domain:
\[
 D(\mathcal{A} ) = \{ u \ | \ u(\bm x) \in L_2(\Omega), \ \sum_{k=0}^{\infty} | (u,\varphi_k) |^2 \lambda_k < \infty \} .
\] 
The operator $\mathcal{A}$ is self-adjoint and positive definite: 
\begin{equation}\label{3}
  \mathcal{A}  = \mathcal{A} ^* \geq \delta I ,
  \quad \delta > 0 ,    
\end{equation} 
where $I$ is the identity operator in $H$. For $\delta$, we have $\delta = \lambda_1$.
In applications, the value of $\lambda_1$ is unknown (the spectral problem must be solved).
Therefore, we assume that $\delta \leq \lambda_1$ in (\ref{3}).
Let us assume for the fractional power of the operator $\mathcal{A}$:
\[
 \mathcal{A}^{\alpha} u =  \sum_{k=0}^{\infty} (u,\varphi_k) \lambda_k^{\alpha }  \varphi_k ,
 \quad 0 < \alpha < 1 . 
\] 
More general and mathematically complete definition of fractional powers of elliptic operators 
is given in \cite{carracedo2001theory,yagi2009abstract}. 

We seek the solution of a Cauchy problem for the evolutionary 
first-order equation with the fractional power of the operator $\mathcal{A}$. 
The solution $u(\bm x,t)$ satisfies the equation
\begin{equation}\label{4}
  \frac{d u}{d t} + \mathcal{A}^{\alpha} u = f(t),
  \quad 0 < t \leq T,  
\end{equation} 
and the initial condition
\begin{equation}\label{5}
  u(0) = u_0 .
\end{equation} 
The key issue in the study of computational algorithms for solving the Cauchy problem (\ref{4}), (\ref{5}) 
is to establish the stability of the approximate solution with respect to small perturbations of the initial 
data and the right-hand side in various norms.

To solve numerically the problem (\ref{4}), (\ref{5}), we employ finite element 
approximations in space (see, e.g., \cite{brenner2008mathematical,Thomee2006}). 
For (\ref{1}) and (\ref{2}), we define the bilinear form
\[
 a(u,v) = \int_{\Omega } \left ( k \, {\rm grad} \, u \, {\rm grad} \, v + c \, u v \right )  d {\bm x} +
 \int_{\partial \Omega } g \, u v d {\bm x} .
\] 
By (\ref{3}), we have
\[
a(u,u) \geq \delta \|u\|^2 .  
\]
Define the subspace of finite elements $V^h \subset H^1(\Omega)$ and the discrete elliptic operator $A$ as
\[
(A y, v) = a(y,v),
\quad \forall \ y,v \in V^h . 
\]
The fractional power of the operator $A$ is defined similarly to $\mathcal{A}^{\alpha}$.
For the spectral problem
\[
A \widetilde{\varphi}_k = \widetilde{\lambda}_k 
\] 
we have
\[
\widetilde{\lambda}_1 \leq \widetilde{\lambda}_2 \leq ... \leq  \widetilde{\lambda}_{M_h},
\quad \| \widetilde{\varphi}_k\| = 1,
\quad k = 1,2, ..., M_h . 
\]
The domain of definition for the operator $A$ is
\[
D(A) = \{ y \ | \ y \in V^h, \ \sum_{k=0}^{M_h} | (y,\widetilde{\varphi}_k) |^2 \widetilde{\lambda}_k < \infty \} .
\]
The operator $A$ acts on the finite dimensional space $V^h$ defined on the domain $D(A)$
and, similarly to (\ref{3}), 
\begin{equation}\label{6}
A = A^* \geq \underline{\delta}_h I ,
\quad \underline{\delta}_h > 0 , 
\end{equation} 
where $\underline{\delta}_h \leq \lambda_1 \leq \widetilde{\lambda}_1$. 
For the fractional power of the operator $A$, we suppose
\[
A^{\alpha } y = \sum_{k=1}^{M_h} (y, \widetilde{\varphi}_k) \widetilde{\lambda}_k^{\alpha} 
\widetilde{\varphi}_k .
\]
The use of finite element approximations for fractional power elliptic operators is discussed in detail, for instance, in
the works \cite{acosta,szekeres2016finite}. 

For the problem (\ref{4}), (\ref{5}), we put into the correspondence the operator equation for $w(t) \in V^h$:
\begin{equation}\label{7}
 \frac{d w}{d t} + A^{\alpha} w = \psi(t), 
 \quad 0 < t \leq T, 
\end{equation} 
\begin{equation}\label{8}
 w(0) = w_0, 
\end{equation} 
where $\psi(t) = P f(t)$, $w_0 = P u_0$ with $P$ denoting $L_2$-projection onto $V^h$.

Now we obtain an elementary a priori estimate for the solution of (\ref{2}), (\ref{3})
assuming that the solution of the problem, coefficients of the elliptic operator, 
the right-hand side and initial conditions are sufficiently smooth.

Let us multiply equation (\ref{2}) by $w$ and integrate it over the domain $\Omega$:
\[
 \left (\frac{d w}{d t}, w \right ) + (A^{\alpha} w, w ) = (\psi, w) .
\]
In view of the self-adjointness and positive definiteness of the operator $A^{\alpha}$, 
we have
\[
 \left (\frac{d w}{d t}, w \right ) \leq  (\psi, w) .
\]
The right-hand side can be evaluated by the inequality
\[
 (\psi, w ) \leq  \|\psi\| \| w \|.
\] 
By virtue of this, we have
\[
 \frac{d}{d t} \|w\| \leq \|\psi \| .
\] 
The latter inequality leads us to the desired a priori estimate:
\begin{equation}\label{9}
 \|w(t)\| \leq \|w_0\| + \int_{0}^{t}\|\psi(\theta) \| d \theta .
\end{equation} 
We will focus on the estimate (\ref{9})
for the stability of the solution with respect to the initial data and the right-hand side 
in constructing discrete analogs of the problem (\ref{7}), (\ref{8}).

\section{Explicit scheme}

To solve numerically the problem (\ref{7}), (\ref{8}), we use simplest explicit and implicit two-level schemes.
Let $\tau$ be a step of a uniform grid in time such that $w^n = w(t^n), \ t^n = n \tau$,
$n = 0,1, ..., N, \ N\tau = T$.
It seems reasonable to begin with the simplest explicit scheme
\begin{equation}\label{10}
 \frac{w^{n+1} - w^{n}}{\tau } + A^\alpha w^{n} = \psi^{n},
 \quad n = 0,1, ..., N-1,
\end{equation} 
\begin{equation}\label{11}
 w^0 = w_0 .
\end{equation}
Advantages and disadvantages of explicit schemes for the standard parabolic problem ($\alpha = 1$)
are well-known, i.e., these are a simple computational implementation and a time step restriction
(see, e.g., \cite{Samarskii1989,SamarskiiMatusVabischevich2002}).
In our case ($\alpha \neq 1$), the main drawback (conditional stability) 
remains, whereas the advantage in terms of implementation simplicity does not exist. 
The approximate solution at a new time-level is determined via (\ref{10}) as
\begin{equation}\label{12}
 w^{n+1} = w^{n} - \tau A^\alpha w^{n} + \tau \psi^{n} . 
\end{equation} 
Thus, we must calculate $A^\alpha w^{n}$.
In view of these problems, considering the scheme (\ref{10}), it is more correct to speak about
the scheme with the explicit approximations in time in contrast to the standard fully explicit scheme.

The numerical implementation of (\ref{12}) is based on the  following representation:
\[
 w^{n+1} = - \tau A r^{n} + w^{n} + \tau \psi^{n} ,
 \quad r^n = A^{\alpha-1} w^{n}.
\] 
We construct a numerical algorithm that employ the rational approximation of the operator
\[
 A^{-\beta}, \quad \beta = 1 - \alpha ,
 \quad 0 < \beta < 1 .  
\]
In this case, we solve standard problems that are related to the operator $A$.

We use an approximation for $A^{-\beta}$ based on integral representation
of a self-adjoint and positive definite operator $A$ (see, e.g., \cite{krasnoselskii1976integral,carracedo2001theory}):
\begin{equation}\label{13}
 A^{-\beta } = \frac{\sin(\pi \beta )}{\pi} \int_{0}^{\infty} \theta^{-\beta } (A + \theta  I)^{-1} d \theta ,
 \quad 0 < \beta  < 1 . 
\end{equation} 
The approximation of $A^{-\beta}$ is based on the use of one or another quadrature formulas for
the right-hand side of (\ref{13}). Various possibilities in this direction are discussed in \cite{bonito2015numerical}.
One possibility is special Gauss-Jacobi quadrature formulas considered in \cite{frommer2014efficient,AcetoNovat}.
Just this approximation of the operator $A^{-\beta}$ is used in the present work.

To achieve higher accuracy in approximating the the right-hand side (\ref{13}),
it is natural to focus on the use of Gauss quadrature formulas.
Some possibilities of constructing quadratures directly for half-infinite intervals are investigated,
for example, in the work \cite{gautschi1991quadrature}.
The classical Gauss quadrature formulas can be used via introducing a new variable of integration.

Suppose (see \cite{frommer2014efficient}) that
\[
 \theta = \mu \frac{1-\eta}{1+\eta},
 \quad \mu > 0.
\] 
From (\ref{13}), we have
\begin{equation}\label{14}
  A^{-\beta } = \frac{2 \mu^{1-\beta} \sin(\pi \beta )}{\pi} \int_{-1}^{1} (1-\eta)^{-\beta}(1+\eta)^{\beta-1}
  \big (\mu (1-\eta) I + (1+\eta) A \big )^{-1} d \eta .
\end{equation}
To approximate the right-hand side of (\ref{14}), we apply the Gauss-Jacobi quadrature formula with the weight
$ (1-\eta)^{\tilde{\alpha}}(1+\eta)^{\tilde{\beta}})$ (see \cite{Rabinowitz}):
\[
  \int_{-1}^{1} f(t) (1-\eta)^{\tilde{\alpha} }(1+\eta)^{\tilde{\beta}} d \eta \approx 
  \sum_{m=1}^{M} \omega_m f(\eta_m) ,
  \quad \alpha, \beta > - 1. 
\] 
Here $\eta_1, \eta_2, ..., \eta_M$ are the roots of the Jacobi polynomial $J_M(\eta; \tilde{\alpha},\tilde{\beta})$
of degree $M$.
The weights $\omega_1, \omega_2, ..., \omega_M$ are given by the formula
\[
\begin{split}
  \omega_m = - & \frac{2 M + \tilde{\alpha} +\tilde{\beta} + 2}{M + \tilde{\alpha} + \tilde{\beta} + 1} 
  \frac{\Gamma (M + \tilde{\alpha} + 1) \Gamma (M + \tilde{\beta} + 1)}{\Gamma (M + \tilde{\alpha} + \tilde{\beta} + 1) (M+1)!} \\
  & \frac{2^{\tilde{\alpha} + \tilde{\beta}}}{J_M^{\,'}(\eta_m; \tilde{\alpha},\tilde{\beta}) 
  J_{M+1}(\eta_m; \tilde{\alpha},\tilde{\beta})} > 0 ,
\end{split}
\] 
where $\Gamma$ denotes the gamma function.

For the fractional power of the operator $A$, from (\ref{14}), we get
\begin{equation}\label{15}
  A^{-\beta } \approx R_M(A),
  \quad  R_M(A) = \sum_{m=1}^{M} d_m (c_m I + A)^{-1} ,
\end{equation} 
where
\[
 \tilde{\alpha} = - \beta ,
 \quad \tilde{\beta} = \beta - 1, 
 \quad d_m = \frac{2 \mu^{1-\beta} \sin(\pi \beta )}{\pi} \frac{\omega_m}{1+\eta_m} ,
 \quad c_m =  \mu  \frac{1-\eta_m}{1+\eta_m} .
\]
In view of (\ref{15}), the approximate solution of the problem $r^n = A^{\alpha-1} w^{n}$
is associated with solving $M$ standard problems $r^n_m = (c_m I + A)^{-1} w^{n}, \ m = 1,2, ..., M$.

Instead of (\ref{12}), we employ the scheme
\begin{equation}\label{16}
  w^{n+1} = - \tau A R_M(A) w^{n} +  w^{n} + \tau \psi^{n} ,
  \quad n = 0, 1, ..., N-1 . 
\end{equation}
Let us consider the stability conditions for the scheme (\ref{11}), (\ref{16}).
For a finite-dimensional self-adjoint operator $A$, in addition to the lower bound (\ref{6}), 
the following upper bound holds:
\begin{equation}\label{17}
 A \leq \overline{\delta}_h I, 
\end{equation} 
where $\overline{\delta}_h = \mathcal{O}(h^{-2})$. Thus
\[
 \underline{\delta}_h^{\alpha}  I \leq A^{\alpha } \leq \overline{\delta}_h^{\alpha }  I ,
 \quad 0 < \alpha < 1 . 
\]
Similar estimates we have  also for $ A R_M(A)$:
\begin{equation}\label{18}
\underline{\gamma}_h I \leq A R_M(A) \leq \overline{\gamma}_h  I ,
 \quad 0 < \alpha < 1 ,
\end{equation}  
with some $\underline{\gamma}_h, \overline{\gamma}_h > 0$. 

\begin{thm}\label{t-1}
If 
\begin{equation}\label{19}
 \tau \leq \tau_0 = \frac{2}{\overline{\gamma}_h} ,
\end{equation} 
then the scheme (\ref{11}), (\ref{16}) is stable in $H$ and
the solution satisfies the following estimate: 
\begin{equation}\label{20}
 \|w^{n+1}\| \leq \|w^{0}\| + \tau \sum_{j=0}^{n} \|\psi^{j}\|,
  \quad n = 0, 1, ..., N-1 .
\end{equation}
\end{thm} 
 
\begin{pf}
From (\ref{16}), we directly obtain
\[
  \|w^{n+1}\| \leq \|I - \tau A R_M(A)\| \|w^{n}\| + \tau \|\psi^{n}\| ,
  \quad n = 0, 1, ..., N-1 . 
\] 
By (\ref{6}), (\ref{17}), (\ref{18}), we get
\[
  \|I - \tau A R_M(A)\| \leq \max_{\underline{\delta}_h \leq z \leq  \overline{\delta}_h} |1 - \tau z R_M(z) | \leq 1
\] 
under the restrictions (\ref{19}). In view of this, we have the level-wise estimate
\[
  \|w^{n+1}\| \leq \|w^{n}\| + \tau \|\psi^{n}\| ,
  \quad n = 0, 1, ..., N-1 , 
\]
that results in the estimate (\ref{20}) for the stability of the solution on the right-hand side and 
the initial conditions.
\end{pf}

It should be noted that the estimate (\ref{20}) for the difference scheme (\ref{11}), (\ref{16}) 
is a discrete analog of the a priori estimate (\ref{9}) for the problem (\ref{7}), (\ref{8}).

Special attention (see \cite{frommer2014efficient,AcetoNovat}) should be given to the problem of choosing the parameter $\mu$ in (\ref{14}).
Taking into account the definition of the operator  $A$, we are interested in the best approximation of $A^{-\beta}$
for the smallest (principal) eigenvalue $\widetilde{\lambda}_1$. In \cite{frommer2014efficient}, there is established 
a remarkable fact that $R_M(z)$ corresponds to a Pade-type approximation for the function $z^{-\beta}$ with expansion point $\mu$.
Thus, the optimal choice corresponds to the selection $\mu = \underline{\delta}_h$. In this case, in (\ref{18}), we have 
$\underline{\gamma} = \underline{\delta}_h^{\alpha}$. The computational complexity of finding $\underline{\delta}_h = \widetilde{\lambda}_1$ 
(the principal eigenvalue of a discrete self-adjoint elliptic operator of second order) is not high.
To this end, it is possible to use standard algorithms (see, e.g., \cite{saad}) and the corresponding software 
(see \cite{hernandez2005slepc}).

The function $z R_M(z)$  for $z \geq z_0 > 0$ is a positive and monotonically increasing function.
In view of this, taking into account (\ref{15}), we have
\[
 \overline{\gamma}_h < \lim_{z \rightarrow \infty} z R_M(z) = \overline{\gamma}(M, \alpha) ,
 \quad \alpha = 1 - \beta , 
\]
where
\[
 \overline{\gamma}(M, \alpha) = \sum_{m=1}^{M} d_m .
\]
From theorem \ref{t-1}, it follows that for the explicit scheme (\ref{11}), (\ref{16}) the time-step restrictions do not depend on
discretization in space, but depend on the power $\alpha$ and the number of approximation nodes $M$. 

\section{Implicit scheme}

Unconditionally stable schemes are constructed on the basis of implicit approximations in time.
Here we consider standard two-level schemes with weights (\cite{Samarskii1989,SamarskiiMatusVabischevich2002}).
For a constant weight parameter $\sigma \ (0 < \sigma \leq 1)$, we define
\[
 w^{n+\sigma} = \sigma w^{n+1} + (1-\sigma) w^{n} .
\]
Instead of (\ref{10}), let us consider the implicit scheme
\begin{equation}\label{21}
 \frac{w^{n+1} - w^{n}}{\tau } + A^\alpha w^{n+\sigma} = \psi^{n+\sigma},
 \quad n = 0,1, ..., N-1 .
\end{equation} 
For $\sigma =1/2$, the difference scheme  (\ref{21}) is the symmetric scheme (the Crank-Nicolson scheme).
It approximates the differential problem 
with the second order by $\tau$, whereas for other values of $\sigma$, we have only the first order.

Rewrite the scheme (\ref{21}) in the form
\[
  \frac{w^{n+\sigma} - w^{n}}{\sigma \tau } + A^\alpha w^{n+\sigma} = \psi^{n+\sigma},
 \quad n = 0,1, ..., N-1 .
\]
In view of this, the transition to a new time-level involves the solution of the problem
\[
 (\nu I + A^\alpha) w^{n+\sigma} = \chi^{n}, 
 \quad \nu = \frac{1}{\sigma \tau} .  
\] 
For this problem, we construct the rational approximation of the operator
\[
 (\nu I + A^\alpha)^{-1},
 \quad 0 < \alpha < 1 . 
\] 

The necessary approximation is based on the following integral representation:
\begin{equation}\label{22}
 (\nu I + A^\alpha)^{-1} = \frac{\sin(\pi \alpha)}{\pi} 
 \int_{0}^{\infty} \frac{\theta^{\alpha}}{\theta^{2\alpha} + 2\theta^{\alpha} \nu  \cos(\pi \alpha) + \nu^2} (A + \theta  I)^{-1} d \theta ,
\end{equation} 
taken from the work \cite{carracedo2001theory} (see Proposition 5.3.2).

Using the new variable $\theta$, from (\ref{22}), we arrive at the representation
\begin{equation}\label{23}
\begin{split}
 (\nu I + & A^\alpha)^{-1} = \frac{2 \mu^{1-\alpha} \sin(\pi \alpha )}{\pi} \\
 & \int_{-1}^{1}    (1-\eta)^{-\alpha}(1+\eta)^{\alpha-1} g(\eta; \nu, \alpha) 
  \big (\mu (1-\eta) I + (1+\eta) A \big )^{-1} d \eta ,
\end{split}
\end{equation} 
where
\[
 g^{-1}(\eta; \nu, \alpha) = 
 1 + 2 \nu  \cos(\pi \alpha) \mu^{-\alpha} \left (\frac{1+\eta}{1-\eta} \right )^\alpha  
 + \nu^2\mu^{-2\alpha} \left (\frac{1+\eta}{1-\eta} \right )^{2\alpha} .  
\]
Further, the Gauss quadrature formula is used (see \cite{Gautschi2004}) with the weight function
\[
 (1-\eta)^{-\alpha}(1+\eta)^{\alpha-1} g(\eta_m; \nu, \alpha) .
\] 
From (\ref{23}), we get
\begin{equation}\label{24}
  (\nu I + A^\alpha)^{-1} \approx R_M(A;\nu),
  \quad  R_M(A;\nu) = \sum_{m=1}^{M} d^\nu_m (c_m I + A)^{-1} .
\end{equation} 
Thereby $R_M(A;0) = R_M(A)$.

For $\sigma > 0$, an approximate solution is obtained from
\begin{equation}\label{25}
\begin{split}
 R^{-1}_M(A;\nu) w^{n+\sigma} & = \nu w^{n} + \psi^{n+\sigma}, \\
 \quad w^{n+1} & = \frac{1}{\sigma} (w^{n+\sigma} - (1-\sigma) w^{n}) ,
 \quad n = 0,1, ..., N-1 .
\end{split} 
\end{equation} 

\begin{thm}\label{t-2}
The difference scheme (\ref{11}), (\ref{25}) for $\sigma \geq 0.5$ and
\begin{equation}\label{26}
 R^{-1}_M(A;\nu) \geq \nu I
\end{equation} 
is unconditionally stable in $H$ and its solution satisfies the a priori estimate
\begin{equation}\label{27}
 \|w^{n+1}\| \leq \|w^{0}\| + \tau \sum_{j=0}^{n} \|\psi^{j+\sigma}\|,
  \quad n = 0, 1, ..., N-1 .
\end{equation}
\end{thm} 
 
\begin{pf}
The use of (\ref{25}) means that instead of (\ref{21}) we employ the scheme
\begin{equation}\label{28}
 \frac{w^{n+1} - w^{n}}{\tau } + D w^{n+\sigma} = \psi^{n+\sigma},
 \quad n = 0,1, ..., N-1 ,
\end{equation} 
where, in view of (\ref{26}), we have
\[
 D = D^* \geq 0,
 \quad   R^{-1}_M(A;\nu) = \nu I + D .
\] 
Multiplying (\ref{28}) scalarly in $H$ by $w^{n+\sigma}$, we get
\begin{equation}\label{28}
 (w^{n+1} - w^{n}, w^{n+\sigma}) \leq \tau (\psi^{n+\sigma}, w^{n+\sigma}) . 
\end{equation} 
Taking into account
\[
 w^{n+\sigma} = (2\sigma -1) w^{n+1} + (1-\sigma) (w^{n+1} + w^{n}) ,
\]  
in view of $0.5 \leq \sigma \leq 1$, for the left-hand side of (\ref{28}), we obtain
\[
\begin{split}
 (w^{n+1} - w^{n}, w^{n+\sigma}) & = (2\sigma -1) \|w^{n+1} \|^2 - (2\sigma -1) (w^{n}, w^{n+1}) \\
 & + (1-\sigma) (\|w^{n+1}\|^2 - \|w^{n}\|^2) \\
 & \geq (\|w^{n+1}\| - \|w^{n}\|) \big ((2\sigma -1) \|w^{n+1} \| + (1-\sigma)(\|w^{n+1}\| + \|w^{n}\|) \big ) \\
 & \geq (\|w^{n+1}\| - \|w^{n}\|) \|w^{n+\sigma}\| .
\end{split}
\] 
For the right-hand side of (\ref{28}), we have
\[
 (\psi^{n+\sigma}, w^{n+\sigma}) \leq \|\psi^{n+\sigma}\| \, \|w^{n+\sigma}\| .
\] 
In view of this, from  (\ref{28}), we get the inequality
\[
 \|w^{n+1}\| \leq  \|w^{n}\| + \tau \|\psi^{n+\sigma}\| ,
\] 
which provides the estimate (\ref{27}).
\end{pf} 

\section{Numerical experiments} 

\begin{figure}
  \begin{center}
  \begin{tikzpicture}[scale = 2]
     \path[draw=black, fill=green!5] (0,0) circle (1.);
     \draw [ultra thick, fill=blue!50] (0,0) -- (0,1) arc [radius=1.0, start angle=90, end angle= 0] -- (1,0) -- (0,0);
     \draw [->] (1,0) -- (1.75,0);
     \draw [->] (0,1) -- (0,1.75);
     \draw(-0.1,-0.1) node {$0$};
     \draw(1.1,-0.1) node {$1$};
     \draw(1.74,-0.1) node {$x_1$};
     \draw(-0.1,1.1) node {$1$};
     \draw(-0.1,1.65) node {$x_2$};
     \draw(0.45,0.45) node {$\Omega$};
     \draw(0.6,-0.13) node {$\Gamma_1$};
     \draw(-0.13,0.6) node {$\Gamma_2$};
     \draw(0.8,0.8) node {$\Gamma_3$};
  \end{tikzpicture}
    \caption{Computational domain $\Omega$}
	\label{f-1}
  \end{center}
\end{figure}
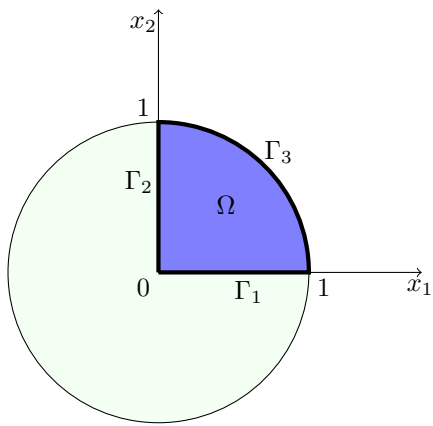 

The test problem under the consideration is constructed using the exact solution of the problem 
in the unit circle (see \cite{vabishchevich2016bnumerical}). 
The computational domain is a quarter of the circle (see Fig.~\ref{f-1}).
Consider the equation
\[
 \mathcal{A} u = - \Delta u ,
 \quad \bm x \in \Omega , 
\]
with the boundary conditions 
\[
 \frac{\partial u}{\partial n} = 0, 
 \quad \bm x \in \Gamma_1,
 \quad \bm x \in \Gamma_2, 
\] 
\[
 \frac{\partial u}{\partial n} + g u = 0,
 \quad \bm x \in \Gamma_3,
 \quad g = \mathrm{const} .  
\]
We study the case, where the solution depends only on 
$r$, and $r = (x_1^2 + x_2^2)^{1/2}$. By virtue of this
\[
 \mathcal{A} u = - \frac{1}{r} \frac{d}{d r} \left ( r  \frac{d u}{d r} \right )  ,
 \quad 0 < r < 1 ,
\] 
\[
 \frac{d u}{d r} + g u = 0,
 \quad r = 1 . 
\] 

The solution of the spectral problem
\[
 - \frac{1}{r} \frac{d}{d r} \left ( r  \frac{d \varphi_k}{d r} \right )  = \lambda_k \varphi_k,
 \quad 0 < r < 1 ,
\] 
\[
 \frac{d \varphi_k}{d r} + g \varphi_k = 0,
 \quad r = 1 , 
\]
is well-known (see, e.g., \cite{Polyanin,Carslaw}).
Eigenfunctions are represented as zero-order Bessel functions:
\[
 \varphi_k(r) = J_0(\sqrt{\lambda_k} r),
\] 
whereas eigenvalues $\lambda_k = \nu_k^2$, where $ \nu_k, \ k = 1,2, ...$ are roots of the equation
\begin{equation}\label{30}
 \nu  J^{\, '}_0(\nu ) + \mu J_0(\nu ) = 0 . 
\end{equation}
The general solution of the homogeneous ($f(t) = 0$) Cauchy problem for  equation (\ref{4}) is
\[
 u(r,t) = \sum_{k=1}^{\infty} a_k \exp(-\nu_k^{2\alpha} t) J_0(\nu_k r) .
\] 

\begin{figure}[htp]
\begin{minipage}{0.33\linewidth}
  \begin{center}
    \includegraphics[width=1.\linewidth] {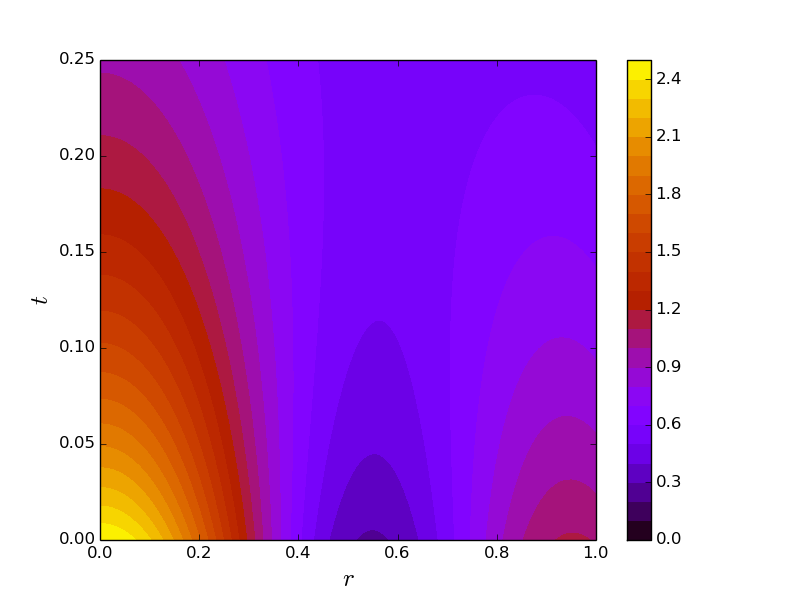}
  \end{center}
\end{minipage}\hfill
\begin{minipage}{0.33\linewidth}
  \begin{center}
    \includegraphics[width=1.\linewidth] {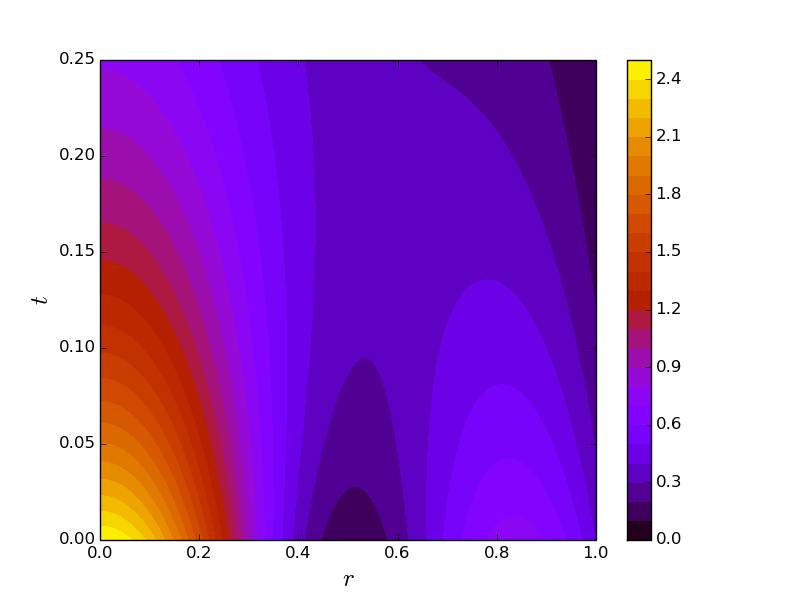}
  \end{center}
\end{minipage}
\begin{minipage}{0.33\linewidth}
  \begin{center}
    \includegraphics[width=1.\linewidth] {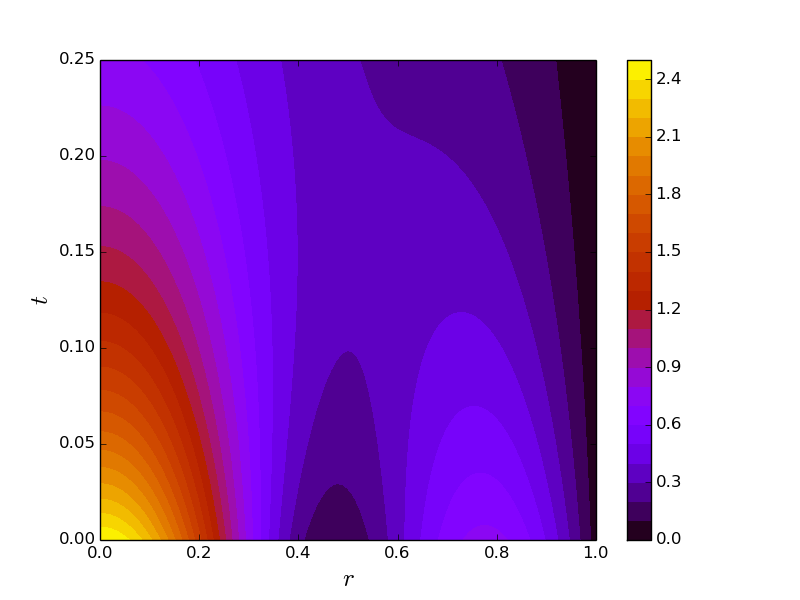}
  \end{center}
\end{minipage}
	\caption{The solution for different $g$ ($\alpha = 0.5$): left: $g = 1$; center:  $g = 10$;
         right: $g = 100$.}
	\label{f-2}
\end{figure}

\begin{figure}[htp]
\begin{minipage}{0.33\linewidth}
  \begin{center}
    \includegraphics[width=1.\linewidth] {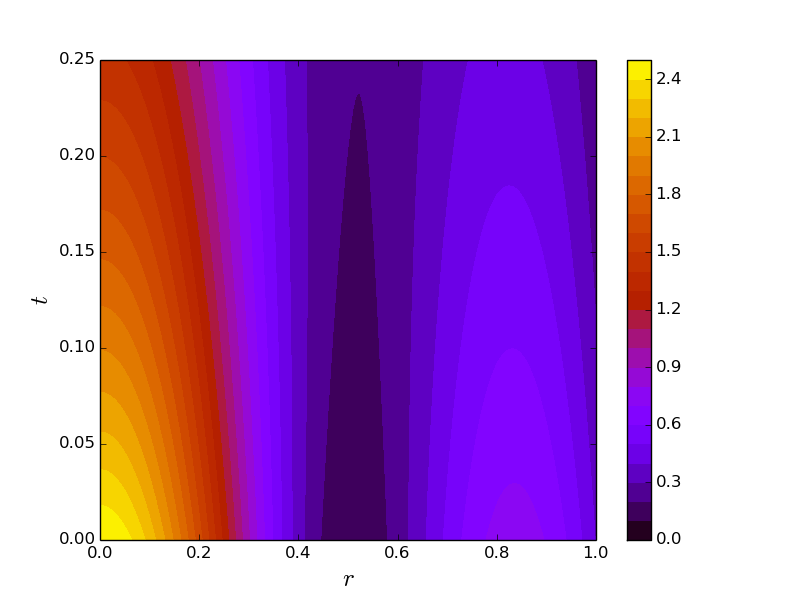}
  \end{center}
\end{minipage}\hfill
\begin{minipage}{0.33\linewidth}
  \begin{center}
    \includegraphics[width=1.\linewidth] {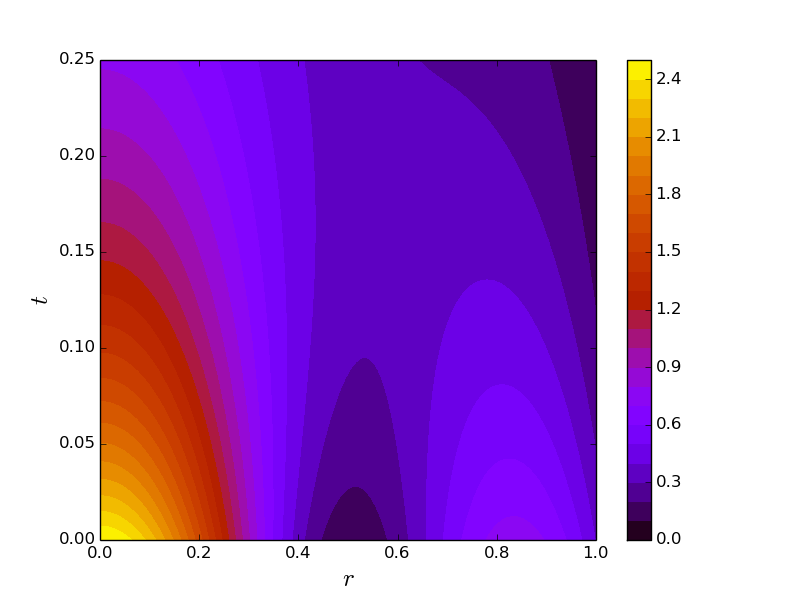}
  \end{center}
\end{minipage}
\begin{minipage}{0.33\linewidth}
  \begin{center}
    \includegraphics[width=1.\linewidth] {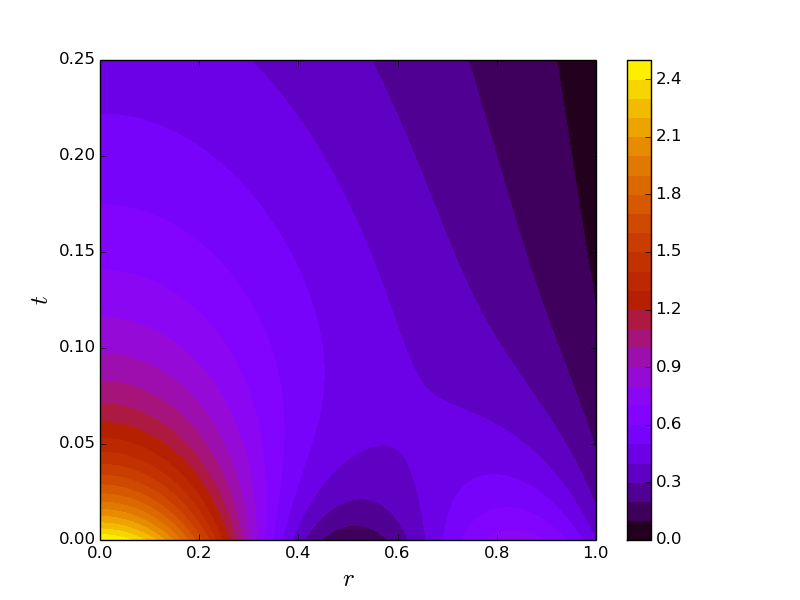}
  \end{center}
\end{minipage}
	\caption{The solution for different $\alpha$ ($g = 10$): left: $\alpha = 0.25$; center:  $\alpha = 0.5$;
         right: $\alpha = 0.75$.}
	\label{f-3}
\end{figure}

To study the accuracy of the approximate solution of the time-dependent equation with the fractional power
of an elliptic operator, we use the exact solution
\begin{equation}\label{31}
 u(r,t) = \exp(-\nu_1^{2\alpha} t) J_0(\nu_1 r) + 1.5 \exp(-\nu_3^{2\alpha} t) J_3(\nu_3 r),
 \quad r = (x_1^2 + x_2^2)^{1/2} .
\end{equation}
The values of the roots $\nu_1, \ \nu_3$ for different values of the boundary condition $\mu$ are given 
in Table~\ref{tab-1}. The exact solution for $T = 0.25$ at different values of $g$ and $\alpha$ is shown in 
Figs.~\ref{f-2} and \ref{f-3}.

\begin{table}
 \caption{The roots of equation (\ref{30})}
 \begin{tabular}{cccc}\label{tab-1}
  $k $ &  $g = 1$  & $g = 10$  & $g = 100$  \\
  \hline
  1 & 1.25578371  & 2.17949660  & 2.38090166 \\
  3 & 7.15579917  & 7.95688342  & 8.56783165 \\
\end{tabular}
\end{table} 

Predictions were performed on a sequence of refining grids, which are shown in Fig.~\ref{4}. 
The numerical solution is compared with the exact one at the final time moment $u(\bm x,T)$. 
Error estimation is performed in $L_2(\Omega)$ and $L_\infty (\Omega)$:
\[
 \varepsilon_2 = \|w^N(\bm x) - u(\bm x,T)\|,
 \quad  \varepsilon_\infty  = \max_{\bm x \in \Omega} |w^N(\bm x) - u(\bm x,T) | .
\] 

\begin{figure}[htp]
\begin{minipage}{0.33\linewidth}
  \begin{center}
    \includegraphics[width=1.\linewidth] {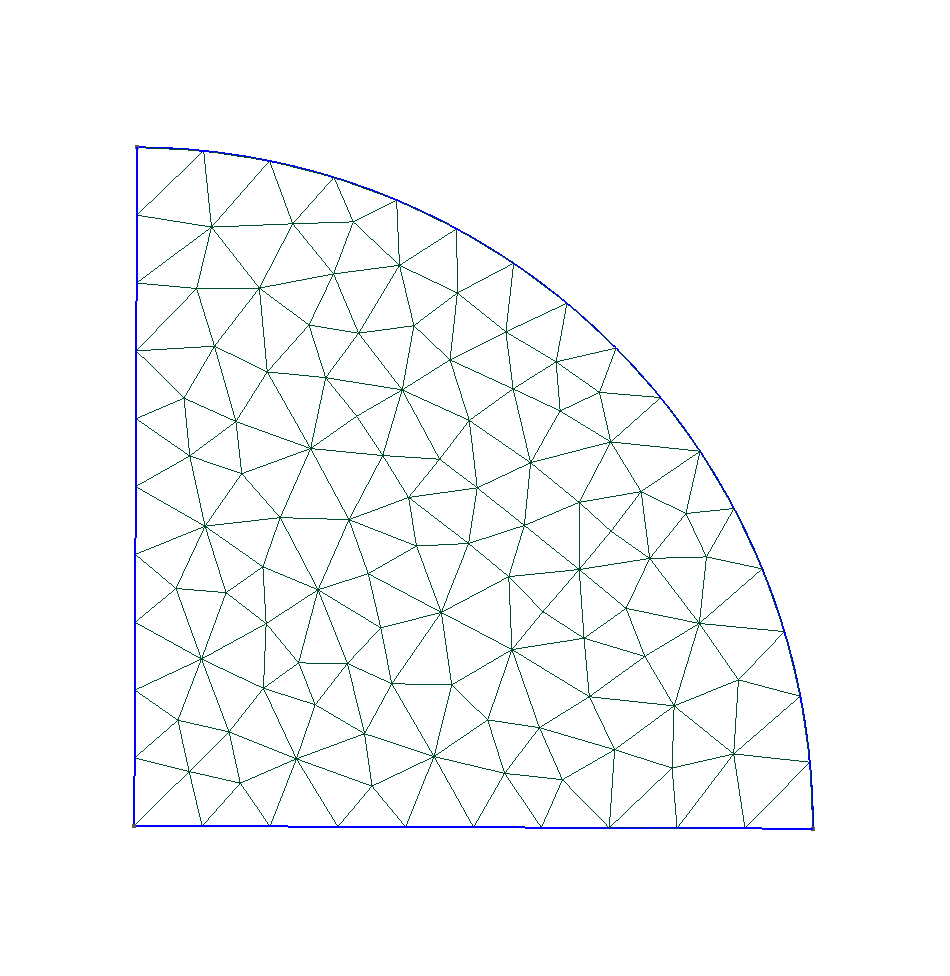}
  \end{center}
\end{minipage}\hfill
\begin{minipage}{0.33\linewidth}
  \begin{center}
    \includegraphics[width=1.\linewidth] {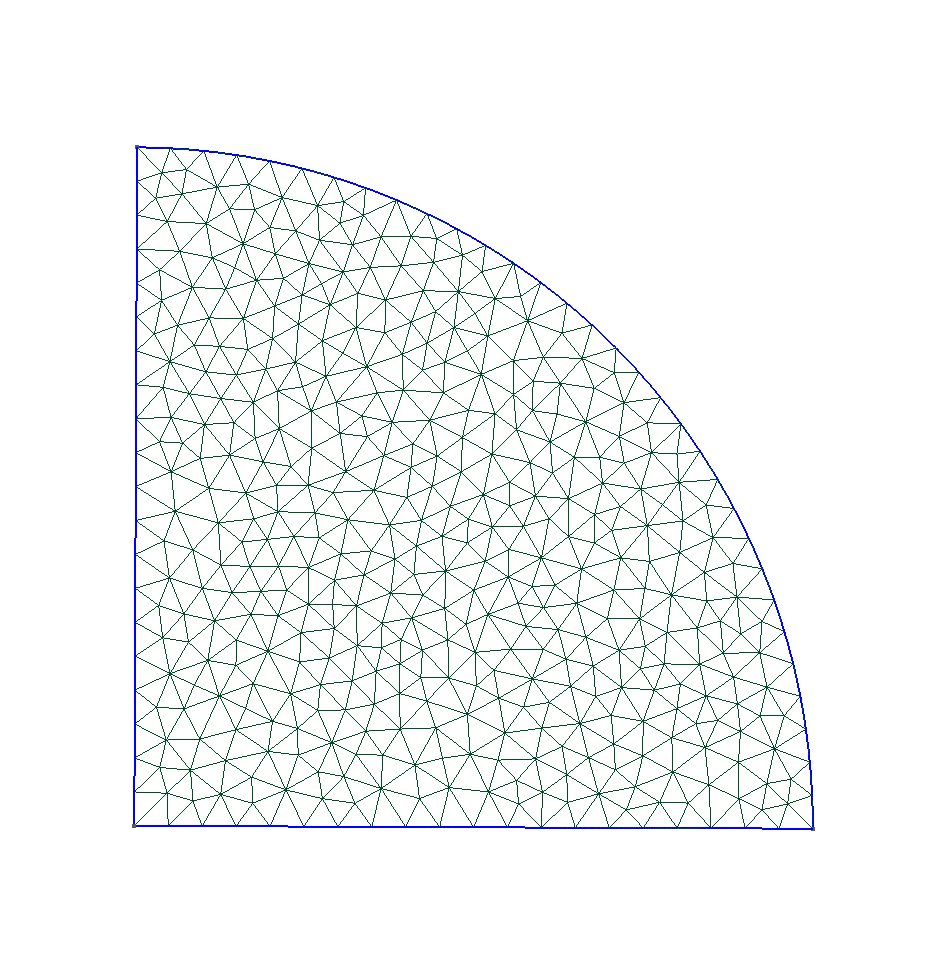}
  \end{center}
\end{minipage}
\begin{minipage}{0.33\linewidth}
  \begin{center}
    \includegraphics[width=1.\linewidth] {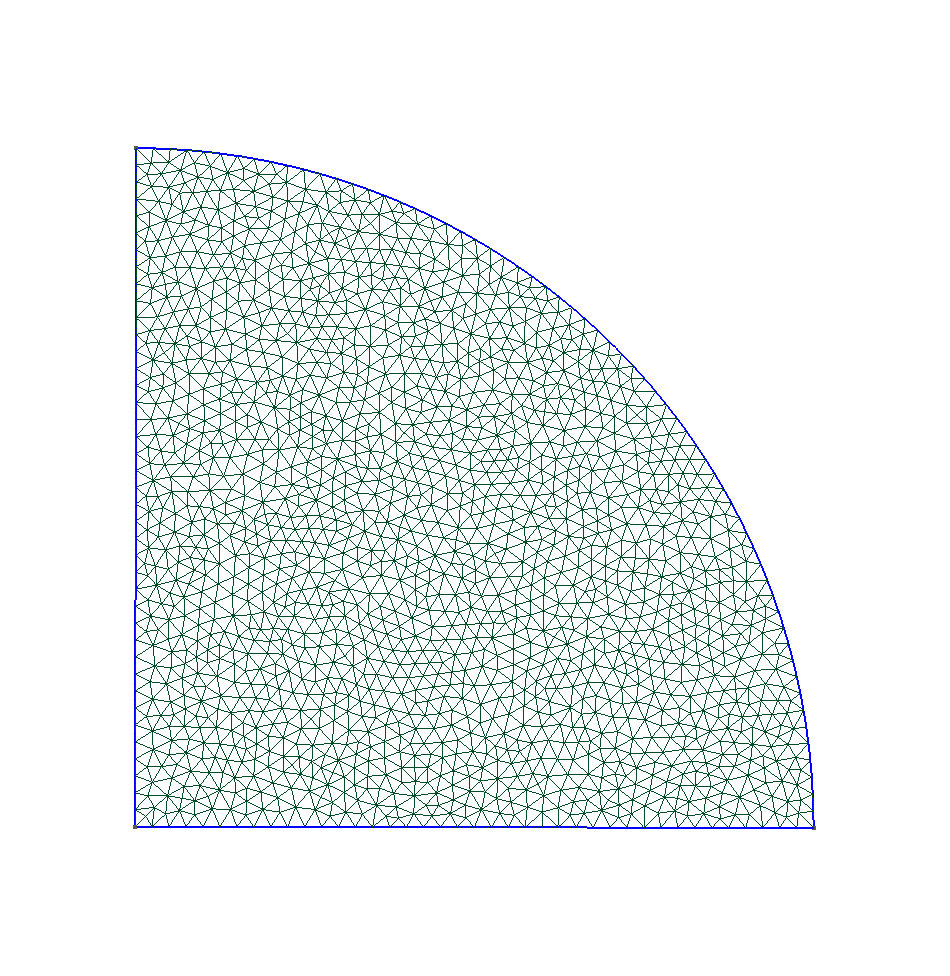}
  \end{center}
\end{minipage}
	\caption{Grid: left: 1 --- 123 vertices and 208 cells; center: 2 --- 461 vertices and 848 cells;
         right: 3 --- 1731 vertices and 3317 cells.}
	\label{f-4}
\end{figure}

\begin{table}[htp]
\caption{The spectrum bounds of the operator $A$}
\label{tab-2}
  \begin{tabular}{c|c|c|c|c}
   $g$  & $\delta = \lambda_1$  & grid & $\underline{\delta}_h$ & $\overline{\delta}_h$ \\
  \hline
      	&  		& 1 &  1.57959231369 & 4225.51507674  \\
  1   	& 1.57699272630 	& 2 &  1.57763558651 & 17104.1780271  \\
  	& 		& 3 &  1.57715815735 & 74989.7519112  \\
  \hline
      	&  		& 1 &  4.76409956820 & 4252.23867499  \\
  10   	& 4.75020542941 	& 2 &  4.75363764524 & 17143.1279728  \\
  	& 		& 3 &  4.75108440807 & 74989.7519112  \\
  \hline
      	&  		& 1 &  5.68846224707 & 7310.80621520  \\
  100   	& 5.66869271459 	& 2 &  5.67358161306 & 22017.7463507  \\
  	& 		& 3 &  5.66994292109 & 74989.7519112  \\
  \end{tabular}
\end{table}

\begin{figure}[htp]
  \begin{center}
    \includegraphics[scale = 0.4] {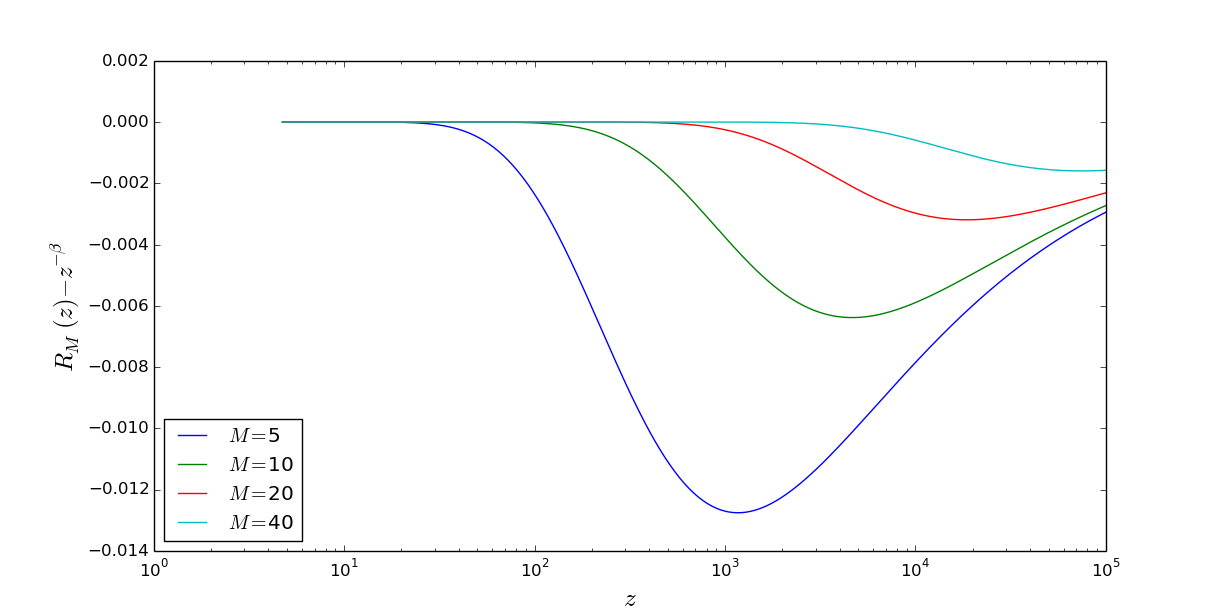}
	\caption{Approximation error for $\beta = 0.5$ ($\mu = \delta$, $z_0^{-\beta}=0.458821546223$).}
	\label{f-5}
  \end{center}
\end{figure} 

The finite element approximation in space is based on the use of continuous $P_1$  Lagrange element, namely,
piecewise-linear elements. The calculations were performed using the computing platform FEniCS 
for solving partial differential equations (website http://fenicsproject.org,  \cite{LoggMardalEtAl2012a,AlnaesBlechta2015a}).
To solve spectral problems with symmetrical matrices, we use the SLEPc library
(Scalable Library for Eigenvalue Problem Computations, http://slepc.upv.es, \cite{hernandez2005slepc}).
We apply the Krylov-Schur algorithm, a variation of the Arnoldi method, proposed by \cite{stewart2002krylov}.

Table~\ref{tab-2} presents the lower and upper bounds of the operator spectrum (see (\ref{6}), (\ref{17}))
on various grids for different values of the parameter $g$ in the boundary condition.
A comparison with the exact minimum eigenvalue demonstrates good accuracy in
evaluation of $\underline{\delta}_h$, which increases with the grid refinement.
It is easy to see a significant dependence of the maximum eigenvalue on the grid size.

Peculiarities of the approximation (\ref{15}) are illustrated by the accuracy of
the approximation of the function $z^{-\beta}$ for $z \geq z_0 = \delta$ (see (\ref{3})).
Figure~\ref{f-5} shows the absolute error arising from the approximation of $z^{-\beta}$ by the function $R_M(z)$
for $\beta = 0.5$ and $g = 10$. In this case, $\mu = \delta$ (see Table~\ref{tab-1}) 
and $R_M(z_0) = z_0^{-\beta}$.  We see higher accuracy near the left boundary $z= z_0$,
whereas for large $z$, the approximation accuracy decreases.
The effect of increasing accuracy with increasing number of nodes of the quadrature formula is clearly observed.

\begin{figure}[htp]
  \begin{center}
    \includegraphics[scale = 0.4] {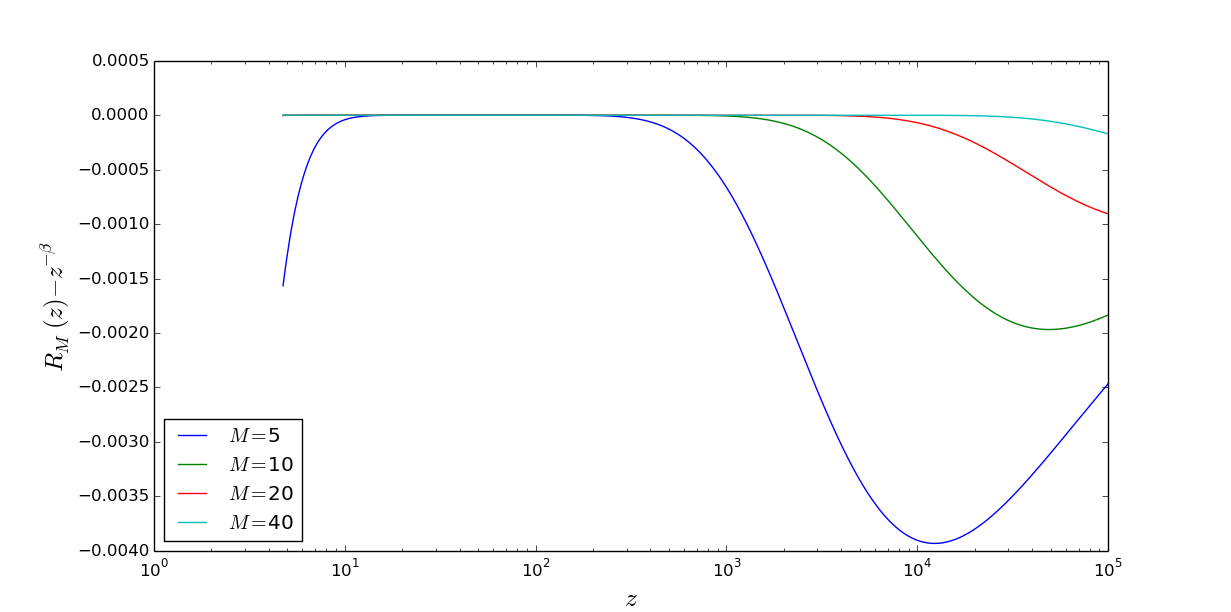}
	\caption{Approximation error for $\beta = 0.5$ ($\mu = 50$).}
	\label{f-6}
  \end{center}
\end{figure}

\begin{figure}[htp]
  \begin{center}
    \includegraphics[scale = 0.4] {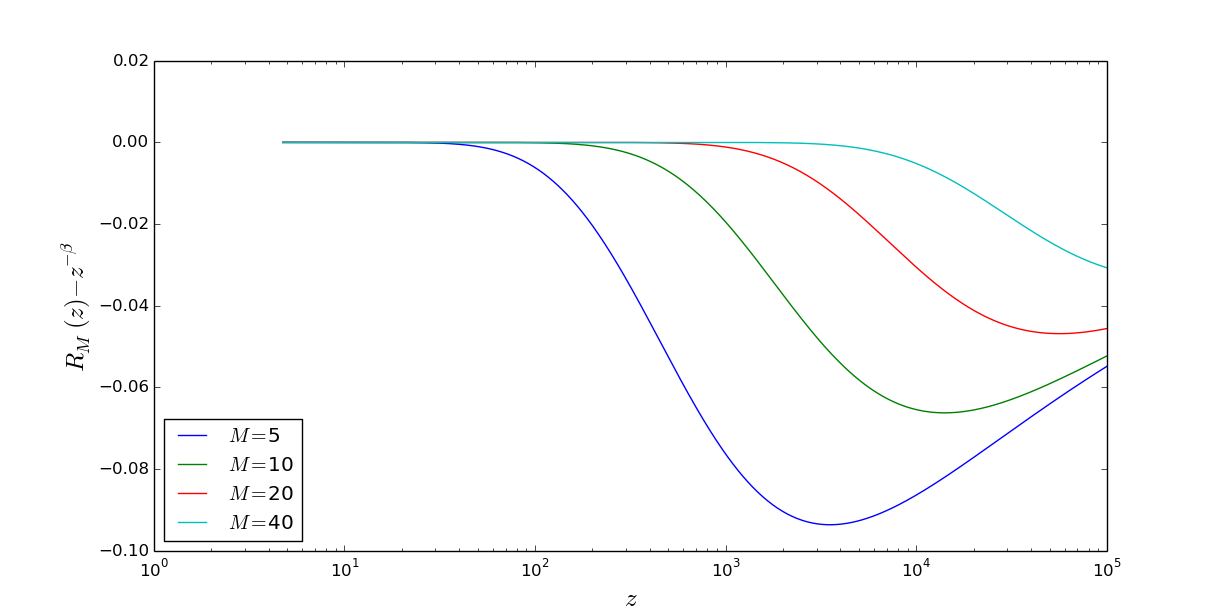}
	\caption{Approximation error for $\beta = 0.25$ ($\mu = \delta$, $z_0^{-\beta}=0.677363673534$).}
	\label{f-7}
  \end{center}
\end{figure}

\begin{figure}[htp]
  \begin{center}
    \includegraphics[scale = 0.4] {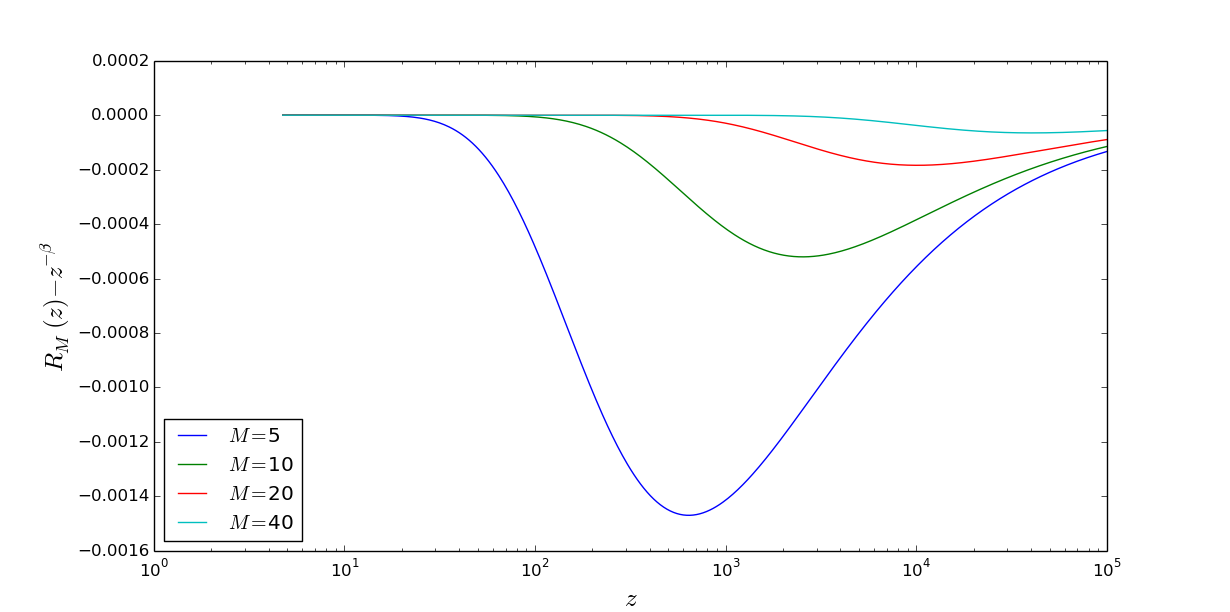}
	\caption{Approximation error for $\beta = 0.75$ ($\mu = \delta$,  $z_0^{-\beta}=0.310789048046$).}
	\label{f-8}
  \end{center}
\end{figure}

Decreasing the approximation accuracy at $z \approx z_0$, we can increase the accuracy for other values of $z$.
For example, Fig.~\ref{f-6} demonstrates the approximation accuracy for $\mu = 50$.
In this case $R_M(\mu) = \mu^{-\beta}$. We need to have good accuracy
for small $z$, and therefore, in calculations, we are guided by the choice of $\mu = \delta$.
The dependence of the approximation accuracy of the function $z^{-\beta}$ on the value of $\beta$ is shown
in Figs.~\ref{f-7} and \ref{f-8}. In these figures, we observe a significant drop in accuracy
with decreasing $\beta$.

\begin{figure}[htp]
  \begin{center}
    \includegraphics[scale = 0.4] {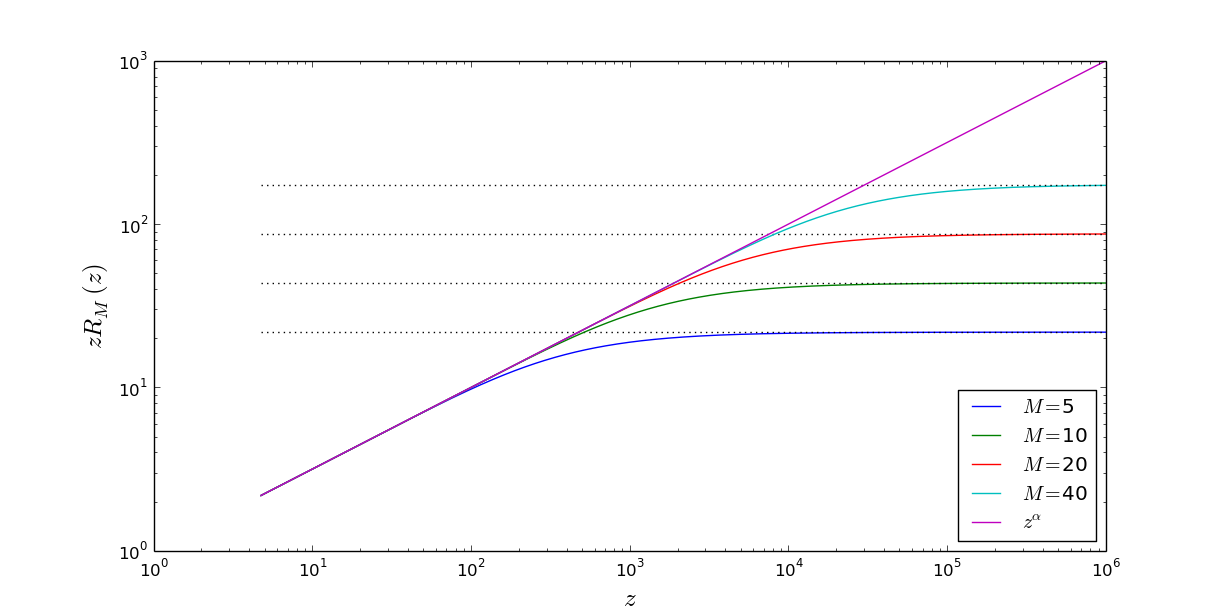}
	\caption{Approximation of $A^\alpha$ ($\alpha = 0.5$) for various $M$.}
	\label{f-9}
  \end{center}
\end{figure}

The numerical implementation of the explicit scheme (\ref{11}), (\ref{16})
involves the approximation of the operator $A^\alpha$ by the expression $A R_M(A) (\beta = 1-\alpha)$.
Peculiarities of this approximation at $\alpha = 0.5$, $g=10$ are shown in Fig.~\ref{f-9}.
It should be noted that the operator $A R_M(A)$ is bounded and the constant $\overline{\gamma}(M,\alpha)$ 
on the right-hand side of (\ref{18})
for the corresponding values of $M$ is presented via dotted lines.

The upper bounds of the operator $A R_M(A)$ are given in Table~\ref{tab-3} for $g=10$.
Increasing $\overline{\gamma}(M,\alpha)$ with increasing the number of quadrature formula nodes $M$ 
results from increasing the accuracy of approximation of the unbounded operator $A^\alpha$.
As $\alpha$ decreases, the value of $\overline{\gamma}(M,\alpha)$ decreases drastically.

\begin{table}[htp]
\caption{The upper bounds of the operator $A R_M(A)$}
\label{tab-3}
  \begin{tabular}{c|c|c|c}
  $M$  & $\overline{\gamma}(M,0.25)$  & $\overline{\gamma}(M,0.5)$ & $\overline{\gamma}(M,0.75)$\\
  \hline
  5	&  4.4602175 &  21.794966  &  142.00220 \\
  10	&  6.3106349 &  43.589932  &  401.45610 \\
  20	&  8.9256294 &  87.179864  &  1135.3565 \\
  40	&  12.623116 &  174.35973  &  3211.1792 \\
  \end{tabular}
\end{table} 

Now we present numerical results obtained using the explicit scheme (\ref{11}), (\ref{16}).
We confine ourselves to the case $\alpha = 0.5$ with the value of the boundary condition parameter $g = 10$.
It is interesting to identify the dependence of accuracy on grids in space and time.
In our case, we should also study the influence of the number of quadrature formula nodes $M$.
Table~\ref{tab-4} demonstrates the numerical solution convergence for decreasing the time step and
increasing the accuracy of approximation of the fractional power operator.
Increasing the accuracy of an approximate solution due to spatial grid refinement is shown in Table~\ref{tab-5}.

\begin{table}[!h]
\caption{Error of the solution for the explicit scheme on grid 2 ($\mu =10, \alpha =0.5$)}
\begin{center}
\begin{tabular}{c|c|c|c|c|c} 
$M$  & $N$ & 25 & 50 & 100 & 200  \\ \hline
5    & $\varepsilon_2$       &  0.00436648  &  0.00207328 &  0.00094905  &  0.00041937 \\
     & $\varepsilon_\infty$  &  0.01896717  &  0.00927482 &  0.00447550  &  0.00216564 \\   \hline
10   & $\varepsilon_2$       &  0.00507635  &  0.00277981 &  0.00164515  &  0.00108352 \\
     & $\varepsilon_\infty$  &  0.02186657  &  0.01217982 &  0.00738292  &  0.00499609 \\   \hline
20   & $\varepsilon_2$       &  0.00507902  &  0.00278251 &  0.00164787  &  0.00108627 \\
     & $\varepsilon_\infty$  &  0.02187724  &  0.01219044 &  0.00739354  &  0.00500671 \\   \hline
40   & $\varepsilon_2$       &  0.00507902  &  0.00278251 &  0.00164787  &  0.00108627 \\
     & $\varepsilon_\infty$  &  0.02187723  &  0.01219043 &  0.00739352  &  0.00500669 \\   
\end{tabular}
\end{center}
\label{tab-4}
\end{table}

\begin{table}[!h]
\caption{Error of the solution for various spatial grids ($\mu =10, \alpha = 0.5, M = 20$)}
\begin{center}
\begin{tabular}{c|c|c|c|c|c} 
grid  & $N$ & 25 & 50 & 100 & 200  \\ \hline
1    & $\varepsilon_2$       &  0.00641387  &  0.00419465 &  0.00310833  &  0.00257568 \\
     & $\varepsilon_\infty$  &  0.02505950  &  0.01634587 &  0.01202969  &  0.00988175 \\   \hline
2    & $\varepsilon_2$       &  0.00507902  &  0.00278251 &  0.00164787  &  0.00108627 \\
     & $\varepsilon_\infty$  &  0.02187724  &  0.01219044 &  0.00739354  &  0.00500671 \\   \hline
3    & $\varepsilon_2$       &  0.00472777  &  0.00241442 &  0.00126921  &  0.00069981 \\
     & $\varepsilon_\infty$  &  0.02077071  &  0.01081355 &  0.00588308  &  0.00342987 \\   
\end{tabular}
\end{center}
\label{tab-5}
\end{table}

\begin{figure}[htp]
  \begin{center}
    \includegraphics[scale = 0.4] {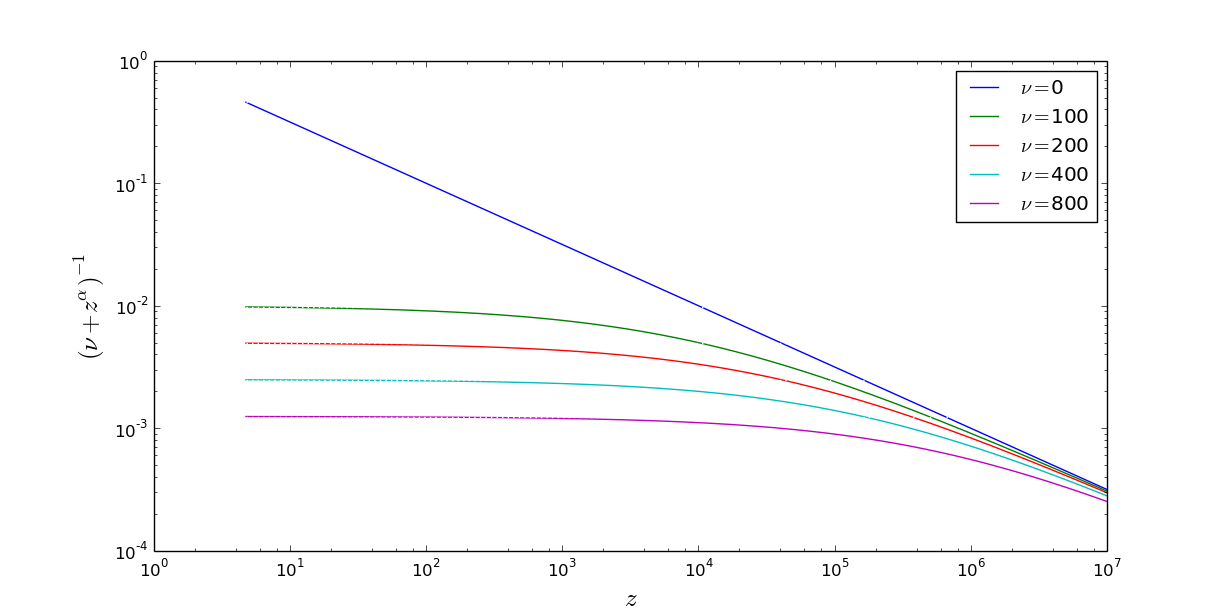}
	\caption{Function $(\nu + z^{\alpha})^{-1}$ for $\alpha = 0.5, z \geq z_0 = \delta$ at various $\nu$.}
	\label{f-10}
  \end{center}
\end{figure} 

The numerical implementation of implicit schemes is associated (see (\ref{21})) with the function
$(\nu + z^{\alpha})^{-1}$. It should be noted that $\nu = (\sigma \tau)^{-1}$ and therefore,
this parameter is large enough in numerical solving unsteady problems.
The function $(\nu + z^{\alpha})^{-1}$, which corresponds to our test problem for $\alpha = 0.5$,
is shown in Fig.~\ref{f-10}. As we noted earlier, if $\nu = 0$, then the optimal value is $\mu = \delta$.
This value is also used in our calculations for $\nu > 0$.

\begin{figure}[htp]
  \begin{center}
    \includegraphics[scale = 0.4] {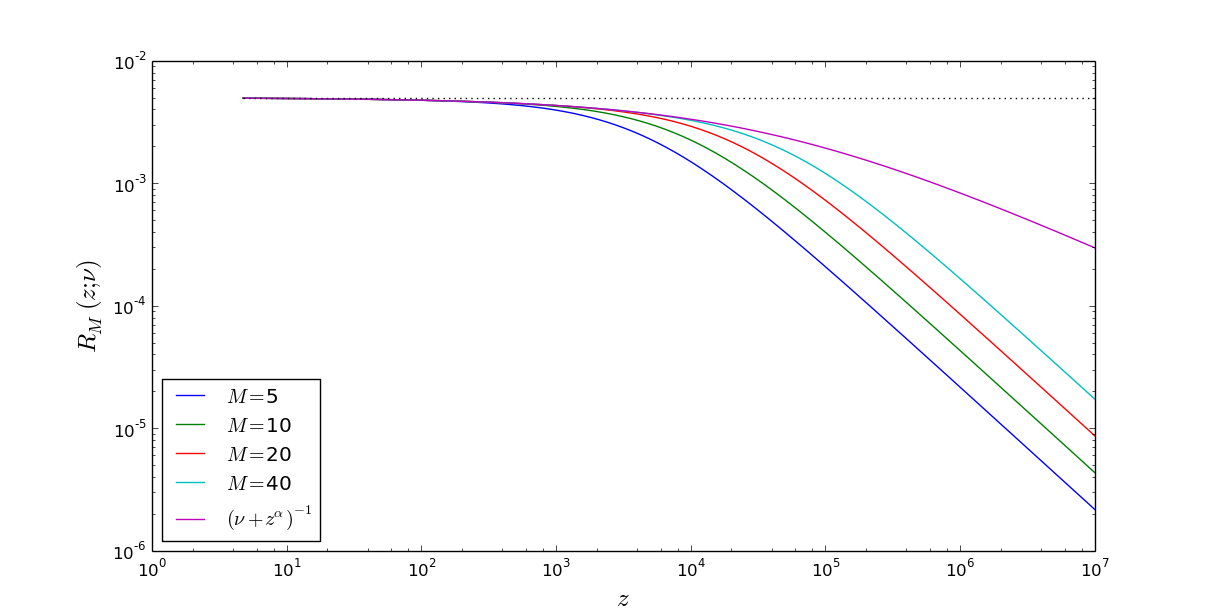}
	\caption{Function $R_M(z;\nu)$ for $\alpha = 0.5, \nu = 200$ at various $M$.}
	\label{f-11}
  \end{center}
\end{figure} 

\begin{figure}[htp]
  \begin{center}
    \includegraphics[scale = 0.4] {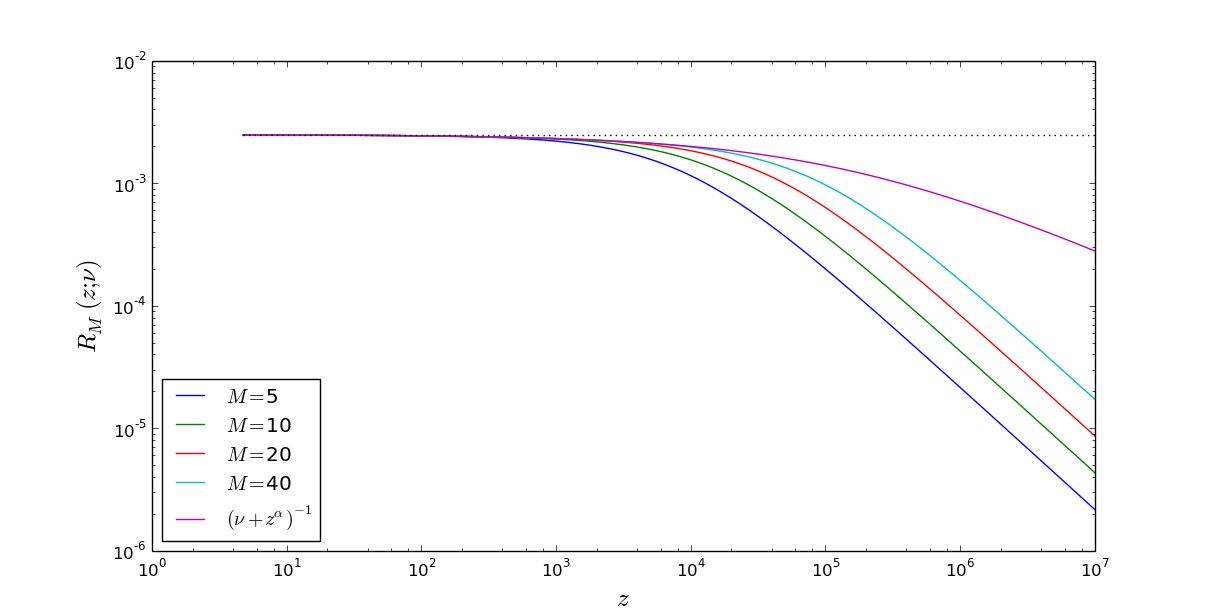}
	\caption{Function $R_M(z;\nu)$ for $\alpha = 0.5, \nu = 400$ at various $M$.}
	\label{f-12}
  \end{center}
\end{figure} 

\begin{figure}[htp]
  \begin{center}
    \includegraphics[scale = 0.4] {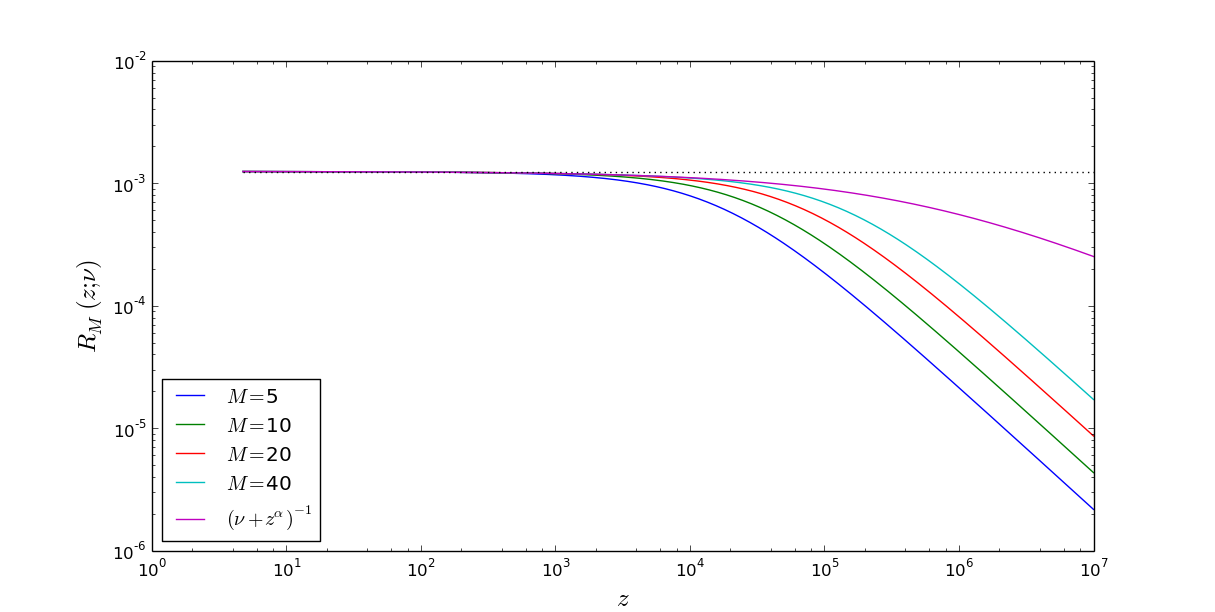}
	\caption{Function $R_M(z;\nu)$ for $\alpha = 0.5, \nu = 800$ at various  $M$.}
	\label{f-13}
  \end{center}
\end{figure} 

Figure~\ref{f-11} shows the approximating function $R_M(z;\nu)$ for $\nu = 200$ with $\mu = \delta$.
The approximation accuracy for various values of $\nu$ is presented in Figs.~\ref{f-12} and \ref{f-13}.
Operator approximations were designed using the package ORTHPOL (see \cite{gautschi1994algorithm}) 
developed for constructing Gauss quadrature formulas with an arbitrary weight function.

The accuracy of the approximate solution of the test problem obtained using the implicit scheme (\ref{11}), (\ref{25}) was investigated
for $\alpha = 0.5$ and $g = 10$.
For the fully implicit scheme ($\sigma = 1$), Table~\ref{tab-6} demonstrates the dependence of the solution accuracy 
on the grid in time for various numbers of the quadrature formula nodes $M$. 

\begin{table}[!h]
\caption{Error of the solution for the implicit scheme on grid 2 ($\mu =10, \alpha =0.5$)}
\begin{center}
\begin{tabular}{c|c|c|c|c|c} 
$M$  & $N$ & 25 & 50 & 100 & 200  \\ \hline
5    & $\varepsilon_2$       &  0.00326281  &  0.00111138 &  0.00044581  &  0.00078676 \\
     & $\varepsilon_\infty$  &  0.01672817  &  0.00699940 &  0.00208775  &  0.00313725 \\   \hline
10   & $\varepsilon_2$       &  0.00394436  &  0.00173870 &  0.00063823  &  0.00018398 \\
     & $\varepsilon_\infty$  &  0.01975034  &  0.01002383 &  0.00511271  &  0.00264674 \\   \hline
10   & $\varepsilon_2$       &  0.00394696  &  0.00174126 &  0.00064056  &  0.00018400 \\
     & $\varepsilon_\infty$  &  0.01976054  &  0.01003442 &  0.00512346  &  0.00265756 \\   \hline
10   & $\varepsilon_2$       &  0.00394696  &  0.00174126 &  0.00064056  &  0.00018399 \\
     & $\varepsilon_\infty$  &  0.01976037  &  0.01003433 &  0.00512339  &  0.00265750 \\   
\end{tabular}
\end{center}
\label{tab-6}
\end{table}

\section*{Conclusion} 

\begin{enumerate}
 \item There is considered a nonclassical problem with the initial data, which is described by
 an evolutionary equation of first order with a fractional power of an elliptic operator. 
 The multidimensional problem is approximated in space using standard 
 finite element piecewise-linear approximations. An a priori estimate for stability with respect to 
 the initial data and the right-hand side is provided.
 \item The explicit scheme is implemented using a Pade-type approximation for the fractional power elliptic operator. 
 Sufficient conditions for the stability of the explicit scheme are formulated. They do not depend on spatial grid steps.
 \item Rational approximation is employed to implement implicit schemes. It is based on a computational generation 
 of Gauss quadrature formulas for an integral representation of the operator of transition to a new time-level.
 \item Possibilities of the proposed algorithms were demonstrated through numerical solving a test two-dimensional problem.
\end{enumerate} 

\section*{Acknowledgements}

The publication was financially supported by the Ministry of Education and Science of 
the Russian Federation (the Agreement \# 02.a03.21.0008).

\end{document}